\RequirePackage{amsmath}
\documentclass[10pt, reqno]{amsart}
\usepackage{amssymb,url, color, mathrsfs}
\usepackage[colorlinks=true, bookmarks=true, pdfstartview=FitH, linktocpage=true, linkcolor = magenta, citecolor = blue]{hyperref}
\usepackage{longtable}
\usepackage{indentfirst, tabularx}
\usepackage[table]{xcolor}
\usepackage{float, hhline, colortbl}
\usepackage{graphicx}

\newtheorem{theorem}{Theorem}
\newtheorem{lemma}{Lemma}

\newtheorem{definition}{Definition}
\newtheorem{remark}{Remark}

\hypersetup{
	colorlinks = true,
	linkcolor = magenta,
	citecolor = blue,
}

\makeatletter
\@addtoreset{equation}{section}
\makeatother

\usepackage{color}

\title[Weighted inequalities]{On some Sobolev and P\'olya-Szeg\"o type inequalities with weights and applications}
\author[T.H.Giang]{Trung Hieu Giang$^{1,2}$}
\address{Trung Hieu Giang$^{1,2}$\newline
	$^1$ Department of Mathematics, City University of Hong Kong, 83 Tat Chee Avenue, Kowloon, Hong Kong; \newline
	$^2$ Institute of Mathematics, Vietnam Academy of Science and Technology, 18 Hoang Quoc Viet, Cau Giay, Hanoi, Vietnam.}
\email{thgiang2-c@my.cityu.edu.hk}

\author[N.M.Tri]{Nguyen Minh Tri$^{2}$}
\address{Nguyen Minh Tri$^{2}$\newline 
	$^2$ Institute of Mathematics, Vietnam Academy of Science and Technology, 18 Hoang Quoc Viet, Cau Giay, Hanoi, Vietnam.}
\email{triminh@math.ac.vn}

\author[D.A.Tuan]{Dang Anh Tuan$^{3}$}
\address{Dang Anh Tuan$^{3}$\newline
	$^3$ University of Sciences, Vietnam National University, 334 Nguyen Trai, Thanh Xuan, Hanoi, Vietnam.}
\email{datuan1105@gmail.com}

\begin{document}

	\begin{abstract}
		We are motivated by studying a boundary-value problem for a class of semilinear degenerate elliptic equations
		\begin{align}\tag{P}\label{P}
			\begin{cases}
				- \Delta_x u - |x|^{2\alpha} \dfrac{\partial^2 u}{\partial y^2} = f(x,y,u) & \textrm{in } \Omega, \\
				u = 0 & \textrm{on } \partial \Omega,
			\end{cases}
		\end{align}
		where $x = (x_1, x_2) \in \mathbb{R}^2$, $\Omega$ is a bounded smooth domain in $\mathbb{R}^3$, $(0,0,0) \in \Omega $, and $\alpha > 0$.
		
		In this paper, we will study this problem by establishing embedding theorems for weighted Sobolev spaces. To this end, we need a new P\'olya-Szeg\"o type inequality, which can be obtained by studying an isoperimetric problem for the corresponding weighted area. Our results then extend the existing ones in \cite{nga, Luyen2} to the three-dimensional context.
	\end{abstract}
	
	\maketitle
	\section{Introduction}\label{sec1}
	
	During the last few decades, degenerate elliptic operators and semilinear degenerate elliptic equations, especially the ones involving Grushin-type operators, have been receiving significant attention from many authors (see, for example, \cite{Luyen, Luyen2, beck, cap, chen, dou}, \cite{Ma}-\cite{Tri3} and the references therein). Among these results, one of the main interests is the question of the existence and nonexistence of nontrivial solutions to the associated (Dirichlet) boundary value problems, with the approach inspired by the classical theory of elliptic partial differential equations. In particular, to establish the nonexistence of nontrivial solutions, the Pohozaev identity is often used as a standard tool. On the other hand, for the existence results, one typically proves corresponding (weighted) Sobolev embeddings and their compactness and invokes the Mountain Pass Lemma. To obtain (weighted) Sobolev embeddings, one can use a P\'olya-Szeg\"o type inequality (the original idea is due to P\'olya and Szeg\"o \cite{PoSz}), which can be established by means of a corresponding (weighted) isoperimetric inequality. One interesting advantage of this approach is that it can give us an (optimal) estimation for the best constant of the Sobolev embeddings. Moreover, it should be noted that weighted inequalities, particularly isoperimetric ones, have also attracted considerable attention (see, for instance, \cite{alvino, aubin, brock, brock2, cab, mon, PoSz} and the references therein).
	
	Recently, based on the aforementioned ideas, Luyen, Tri, and Tuan \cite{Luyen2} studied the following problem
	\begin{align}\label{eqn1}
		\begin{cases}
			- \dfrac{\partial^2u}{\partial x^2} - |x|^{2\alpha} \dfrac{\partial^2 u}{\partial y^2} = f(x,y,u) & \textrm{in } \Omega, \\
			u = 0 & \textrm{on } \partial \Omega,
		\end{cases}
	\end{align}
	where $\Omega$ is a bounded smooth domain in $\mathbb{R}^2$, $(0,0) \in \Omega $, and $\alpha > 0$. Their research focused on the question of whether nontrivial solutions exist or not. One of the key ingredients in their paper is the weighted embedding theorem, which Nga, Tri, and Tuan \cite{nga} established via a P\'olya–Szeg\"o type inequality with respect to weighted areas. However, we would like to emphasize that the results in \cite{nga, Luyen2} are only applicable to the two-dimensional case (see Remark \ref{rem1}), and thus, it is a natural question whether they can be generalized to a higher-dimensional context. 
	
	In this paper, our purpose is to extend the results given in \cite{nga, Luyen2} to the three-dimensional context. More specifically, one objective of this paper, achieved in Sections \ref{sec2} and \ref{sec3}, is to establish weighted inequalities starting with the following isoperimetric inequality.
	
	\begin{theorem}[The weighted isoperimetric inequality]\label{thm1}
		Let $E \subset \mathbb{R}^3$ be a nonempty bounded open set with Lipschitz boundary $\partial E$ such that $P_{2,\alpha} (E)$ is finite. Then the following inequality holds for every $j=1,...,2n(\alpha)$
		\begin{equation}\label{eqn2}
			\dfrac{P_{2,\alpha,j}(B_{s_j})^{3/2}}{|B_{s_j}|_{2,\alpha}} \leq \dfrac{P_{2,\alpha}(E)^{3/2}}{|E|_{2,\alpha}},
		\end{equation}
		where $B_{s_j} := \left\{ (x_1,x_2,y) \in \mathbb{R}^3_{s_j}, \textrm{ }\dfrac{|x|^{2\alpha+2}}{(\alpha+1)^2} + y^2 < 1 \right\}$. Here, the notations $P_{2,\alpha}$, $P_{2,\alpha,j}$, and $|\cdot|_{2,\alpha}$ respectively denote weighted areas and weighted volume, which will be precisely given in Definitions \ref{def1}-\ref{def3}. 
	\end{theorem}
	
	Then, using this, we will be able to prove the following P\'olya-Szeg\"o type inequality.
	
	\begin{theorem}[The P\'olya-Szeg\"o type inequality]\label{thm2}
		Let $u \in C^\infty_0 ( \mathbb{R}^3; \mathbb{R}_+)$. Then
		$$ \int\limits_{\mathbb{R}^3_{s_1}} \left| \nabla_G u^*\right|^2dx_1dx_2dy \leq \int\limits_{\mathbb{R}^3}\left| \nabla_G u\right|^2dx_1dx_2dy,$$
		where $\nabla_G u := \big(\frac{\partial u}{\partial x_1}, \frac{\partial u}{\partial x_2}, |x|^\alpha \frac{\partial u}{\partial y} \big)$ and $\left|\nabla_G u \right| := \left( |\frac{\partial u}{\partial x_1}|^2 + |\frac{\partial u}{\partial x_2}|^2 + |x|^{2\alpha} |\frac{\partial u}{\partial y}|^2 \right)^{1/2}$. Here, $u^*$ denotes the rearrangement of $u$, and will be defined in Definition \ref{def4}.
	\end{theorem}
	
	The above inequality plays a crucial role in the proof of the following weighted Sobolev embedding.
	
	\begin{theorem}\label{thm3}
		The following inequality holds
		\begin{equation}\label{eqn3}
			\Bigg(\int\limits_{\mathbb{R}^3} |x|^{2\alpha}|u|^{6} dx_1dx_2dy \Bigg)^{1/6} \leq C^{-1}_{2\alpha,6}\Bigg(\int\limits_{\mathbb{R}^3}|\nabla_G u|^2 dx_1dx_2dy\Bigg)^{1/2}
		\end{equation}
		for every $u \in W^{1,2\alpha,6}_0(\mathbb{R}^3)$. Here, the defintion of the space $W^{1,2\alpha,6}_0(\mathbb{R}^3)$ will be given in Section \ref{sec3}.
		
		Furthermore, the best constant $C_{2\alpha,6}$ satisfies
		$$C_{2\alpha,6} \geq \left( \dfrac{2\pi}{n(\alpha)}\right)^{\frac{-1}{3}} (\alpha+1)^{\frac{-1}{3}} D_{2,6,3},$$
		where $D_{2,6,3}$ is a constant which will be given in Lemma \ref{lem2}. 
	\end{theorem}
	
	The other objective of this paper, achieved in Section \ref{sec4}, is to establish the existence and nonexistence results for the problem \eqref{P}. In the case where the function $f$ is given by
	\begin{equation}
		\label{eqn4}
		f(x_1,x_2,y,\xi) := |x|^{2\alpha}|\xi|^{p-1}\xi,
	\end{equation}
	the nonexistence of nontrivial solutions will be derived from a Pohozaev-type identity. The statement of this result is as follows.
	
	\begin{theorem}\label{thm4}
		Assuming $f$ has the form given in \eqref{eqn4}. Let $\Omega$ be $G_\alpha$ -- star-shaped with respect to the origin and $p>5$ (the definition of this type of domains will be given in Definition \ref{def6}). Then, the problem \eqref{P} has no nontrivial solution $u \in \mathcal{S}_2 (\overline{\Omega})$. Here, the function space $\mathcal{S}_2 (\overline{\Omega})$ will be defined in the beginning of Subsection \ref{sec41}.
	\end{theorem}
	
	Next, regarding existence results, we need the following assumptions: Assume that $f: \Omega \times \mathbb{R} \to \mathbb{R}$ is a Carath\'eodory function such that
	\begin{itemize}
		\item[(A1)] $f$ satisfies the assumption in Lemma \ref{lem4} and the additional condition 
		\begin{equation}\label{eqn5}
			p_1 > \dfrac{3}{2};
		\end{equation}
		\item[(A2)] there exist $C \in [0, +\infty)$ and $\psi \in L^1_{|x|^{2\alpha}} (\Omega)$ (for the notation $L^1_{|x|^{2\alpha}} (\Omega)$ see Definition \ref{def7}) such that 
		$$|f(x_1,x_2, y, \xi)| \leq |x|^{2\alpha} \psi(x_1,x_2, y)$$ 
		for a.e. $(x_1,x_2, y)$ in $\Omega$ and for every $|\xi| \leq C$;
		\item[(A3)] there exists a non-positive function $\varphi$ such that 
		$$\int_\Omega |\varphi(x_1,x_2, y)|dx_1dx_2dy < \infty$$
		and 
		$$\varphi(x_1,x_2, y) \leq \frac{ f(x_1,x_2,y,\xi)}{ \xi}$$
		for a.e. $(x_1,x_2, y) \in \Omega$ and for every $\xi \in \mathbb{R}_+$;
		\item[(A4)] $f(x_1,x_2, y, 0) = 0$ for a.e. $(x_1,x_2, y)$ in $\Omega$ and the following limits hold uniformly for a.e. $(x_1,x_2, y)$ in $\Omega$
		$$\lim\limits_{\xi\to0} \dfrac{f(x_1,x_2, y, \xi)}{ |x|^{2\alpha} \xi} = 0 \textrm{ and }\lim\limits_{ |\xi|\to+\infty} \dfrac{f(x_1,x_2, y, \xi)}{ |x|^{2\alpha}\xi}  = +\infty;$$
		\item[(A5)] $\frac{f(x_1,x_2,y,\xi )}{\xi}$ is increasing in $\xi \geq C$ and decreasing in $\xi \leq -C$ for a.e. $(x_1,x_2, y)$ in $\Omega$.
	\end{itemize}
	
	Notice that the above assumptions are natural in the context of degenerate elliptic problems, as they generalize the ones usually given in the context of elliptic problems (see, for example, \cite{Liu}). Under these assumptions, we have the following existence result.
	
	\begin{theorem}\label{thm5}
		Suppose that $f$ satisfies (A1)-(A5). Then, the problem \eqref{P} has a nontrivial weak solution (see Definition \ref{def10} for the definition of weak solutions).
	\end{theorem}
	
	The organization of this paper is as follows. In Section \ref{sec2}, we will give the proofs for Theorems \ref{thm1} and \ref{thm2}. Then, in Section \ref{sec3}, we will demonstrate how Theorem \ref{thm2} will be utilized to establish the Sobolev-type inequality with weight presented in Theorem \ref{thm3}. Here, we would like to emphasize that our method, based on \cite{nga}, is only able to give lower estimates for the best constants appearing in the weighted Sobolev inequalities, and thus calculating these constants remains an open question. Finally, in Section \ref{sec4} we will apply the above results and the Pohozaev-type identity to give the proofs for Theorems \ref{thm4} and \ref{thm5}.
	
	\section{A P\'olya-Szeg\"o type inequality}\label{sec2}
	In this section, we will prove our new weighted isoperimetric and P\'olya-Szeg\"o type inequalities mentioned in Section \ref{sec1}. To this end, we first introduce some notations.
	
	Throughout this paper, let $\alpha$ be a fixed given positive real number and $n(\alpha)$ be the smallest positive integer such that $n(\alpha) \geq \alpha+1$. We denote 
	$$\mathbb{R}^3_{s_j} := \left\{(x_1,x_2,y) = (r\cos\theta, r\sin\theta, y) \in \mathbb{R}^3, \textrm{ } r >0, \textrm{ } \theta \in \Big(\frac{(j-1)\pi}{n(\alpha)},\frac{j\pi}{n(\alpha)}\Big)\right\},$$
	for each $j=1,...,2n(\alpha)$. 
	
	Let $E \subset \mathbb{R}^3$ be a bounded open set with Lipschitz boundary $\partial E$. We denote by $\nu = (\nu_1, \nu_2, \nu_3)$ the outward unit normal to $\partial E$.  We also denote
	$$ \partial_{s_j} E:= \partial E \cap \mathbb{R}^3_{s_j}$$
	for each $j=1,...,2n(\alpha)$, and 
	$$ |x| := \sqrt{x_1^2 + x_2^2} \textrm{ for each } (x_1,x_2) \in \mathbb{R}^2.$$
	
	\begin{definition}\label{def1}
		The $(2,\alpha)$-volume of $E$ is defined by
		$$ |E|_{2,\alpha} = \int\limits_E |x|^{2\alpha} dx_1dx_2 dy.$$ 
	\end{definition}
	
	\begin{definition}\label{def2}
		The $(2,\alpha)$-area of $E$ is defined by
		$$ P_{2,\alpha} (E) = \int\limits_{\partial E} |x|^\alpha \sqrt{ \nu_1^2 + \nu_2^2 + |x|^{2\alpha} \nu_3^2}d\mathcal{H}^2,$$ 
		where $\mathcal{H}^2$ denotes the two-dimensional Hausdorff surface measure in $\mathbb{R}^3$.
	\end{definition}
	
	\begin{definition}\label{def3}
		For each $j=1,...,2n(\alpha)$, the $(2,\alpha,j)$-area of $E$ is defined by
		$$ P_{2,\alpha,j} (E) = \int\limits_{\partial_{s_j} E} |x|^\alpha \sqrt{ \nu_1^2 + \nu_2^2 + |x|^{2\alpha} \nu_3^2}d\mathcal{H}^2.$$ 
	\end{definition}
	
	Notice that
	$$P_{2,\alpha}(E_\lambda) = \int\limits_{\partial E_\lambda} |X|^\alpha\sqrt{\mu_1^2 +\mu_2^2 + |X|^{2\alpha} \mu_3^2}d\mathcal{H}^2 = \lambda^{2\alpha+2} P_{2,\alpha} (E), $$
	$$P_{2,\alpha,j}(E_\lambda) = \int\limits_{\partial_{s_j} E_\lambda} |X|^\alpha\sqrt{\mu_1^2 +\mu_2^2 + |X|^{2\alpha} \mu_3^2}d\mathcal{H}^2 = \lambda^{2\alpha+2} P_{2,\alpha,j} (E), $$
	and
	$$ |E_\lambda|_{2,\alpha} = \int\limits_{E_\lambda} |X|^{2\alpha} dX_1dX_2dY = \lambda^{3\alpha+3} |E|_{2,\alpha},$$
	where $E_\lambda = \{ (\lambda x_1, \lambda x_2, \lambda^{\alpha+1} y): (x_1, x_2, y) \in E \}$ and $\mu = (\mu_1, \mu_2, \mu_3)$ is the outward unit normal to $\partial E_\lambda$. Thus, the fractions
	$$\dfrac{P_{2,\alpha}(E)^{3/2}}{|E|_{2,\alpha}} \textrm{ and } \dfrac{P_{2,\alpha,j}(E)^{3/2}}{|E|_{2,\alpha}}$$ 
	for each $j=1,...,n(\alpha)$, are invariant under the scaling $(x_1,x_2,y) \mapsto (\lambda x_1, \lambda x_2, \lambda^{\alpha+1} y)$.
	
	The proof of Theorem \ref{thm1} is given as follows.
	
	\begin{proof}[Proof of Theorem \ref{thm1}]
		For readers' convenience, we will divide the proof into two steps.
		
		\noindent\textbf{Step 1:} Assuming that there exists an integer $j$ such that $1\leq j\leq 2n(\alpha)$ and $E \subset \mathbb{R}^3_{s_j}$ with Lipschitz boundary. We will show that
		\begin{equation}\label{eqn6}
			\dfrac{P_{2,\alpha,j}(B_{s_j})^{3/2}}{|B_{s_j}|_{2,\alpha}} \leq \dfrac{P_{2,\alpha,j}(E)^{3/2}}{|E|_{2,\alpha}}.
		\end{equation}
		
		Without loss of generality, we can assume that $j=1$. We define $\Phi_1: \mathbb{R}_+ \times \left( 0, \dfrac{\pi}{n(\alpha)} \right) \times \mathbb{R} \to \mathbb{R}^3_{s_1}$ by
		$$ \Phi_1 (r, \theta, y) = (r\cos \theta, r\sin \theta, y).$$
		It is easy to see that $\Phi_1$ is a homeomorphism. Next, we define $\Phi_2: \mathbb{R}_+ \times \left( 0, \dfrac{\pi}{n(\alpha)} \right) \times \mathbb{R} \to \mathbb{R}^3_{\widetilde{s}_1}$ by
		$$ \Phi_2 (r, \theta, \eta) = \left(\dfrac{r^{\alpha+1}\cos (\alpha+1)\theta}{\alpha+1}, \dfrac{r^{\alpha+1}\sin(\alpha+1)\theta}{\alpha+1}, \eta\right), $$
		where
		$$\mathbb{R}^3_{\widetilde{s}_1} := \left\{(\xi_1,\xi_2,\eta) = (\rho\cos\varphi, \rho\sin\varphi, y) \in \mathbb{R}^3, \textrm{ } \rho >0, \textrm{ } \varphi \in \Big(0,\frac{(\alpha+1)\pi}{n(\alpha)}\Big) \right\}.$$
		It is also easy to verify that $\Phi_2$ is a homeomorphism. Hence
		$$ \Psi = \Phi_2 \circ \Phi_1^{-1} $$
		is a homeomorphism. 
		
		Next, let $\widetilde{E} = \Psi (E)$ and $\widetilde{B}_{s_1} = \Psi(B_{s_1})$. We deduce that
		$$ \widetilde{B}_{s_1} = \left\{ (\xi_1, \xi_2, \eta) \in \mathbb{R}^3_{\widetilde{s}_1}, \textrm{ } \xi_1^2 + \xi_2^2 + \eta^2 < 1\right\}.$$
		By some straightforward calculations, we obtain
		$$ |E|_{2,\alpha} = \int\limits_{\widetilde{E}} d\xi_1d\xi_2d\eta = |\widetilde{E}|, \quad |B_{s_1}|_{2,\alpha} = \int\limits_{\widetilde{B}_{s_1}} d\xi_1d\xi_2d\eta = |\widetilde{B}_{s_1}|, $$
		$$ P_{2,\alpha,1} (E) = \int\limits_{\partial_{\widetilde{s}_1} \widetilde{E}} d\mathcal{H}^2 =: P_{s_1}(\widetilde{E}), \quad P_{2,\alpha,1} (B_{s_1}) = \int\limits_{\partial_{\widetilde{s}_1} \widetilde{B}_{s_1}} d\mathcal{H}^2 =: P_{s_1}( \widetilde{B}_{s_1}),$$
		where
		$$ \partial_{\widetilde{s}_1} \widetilde{E} = \partial \widetilde{E} \cap \mathbb{R}^3_{\widetilde{s}_1}, \quad \partial_{\widetilde{s}_1} \widetilde{B}_{s_1} = \partial \widetilde{B}_{s_1} \cap \mathbb{R}^3_{\widetilde{s}_1}.$$
		
		Now, by applying \cite[Theorem 1.1]{Lions}, we deduce that
		$$ \dfrac{P_{s_1}(\widetilde{E})^{3/2}}{|\widetilde{E}|} \geq \dfrac{P_{s_1}(\widetilde{B}_{s_1})^{3/2}}{|\widetilde{B}_{s_1}|}, $$
		and thus \eqref{eqn6} follows. 
		
		\noindent \textbf{Step 2:} For the open, bounded set $E \subset \mathbb{R}^3$ with Lipschitz boundary, we set
		$$E = \bigcup_{i=1}^{2n(\alpha)} E_{s_i}, $$
		where $E_{s_i} := E \cap \overline{\mathbb{R}}^3_{s_i}$. It is easy to see that
		\begin{equation}\label{eqn7}
			\bigcup_{i=1}^{2n(\alpha)} \partial_{s_i} E \subset \partial E.
		\end{equation}
		Notice that the sets $E_{s_i}$ are respectively open, bounded subsets in $\overline{\mathbb{R}}^3_{s_i}$ with Lipschitz boundaries. Hence we have the inequality \eqref{eqn6} for each $E_{s_i}$. Since $3/2 > 1$, we have that
		$$ P_{2,\alpha} (E)^{3/2} \geq \sum\limits_{j=1}^{2n(\alpha)} P_{2,\alpha,j}( E_{s_j})^{3/2}.$$
		By combining this and the inequality \eqref{eqn6} for every $E_{s_i}$, our proof is complete. 
	\end{proof}
	
	\begin{remark}
		\label{rem1}
		The idea of our proof is similar to the one given in \cite{nga}, which was based on a change of variables, so that the authors in \cite{nga} were able to apply the result on a weighted isoperimetric inequality given in \cite{cab}. However, when passing from the two-dimensional case to the three-dimensional one, we need to modify the proof, and we are only able to obtain a result for the case of the weight $|x|^{p\alpha}$ with $p=2$, but not for general $p$. The main reason is that for the two-dimensional case, the weight $|x|^{p\alpha}$ is already monomial, which allows us to apply the result in \cite{cab}. For the higher-dimensional case, however, this is no longer true.
	\end{remark}

	\begin{definition}\label{def4}
		Let $u \in C_0^\infty (\mathbb{R}^3; \mathbb{R}_+)$. The rearrangement $u^* : \overline{\mathbb{R}}^3_{s_1} \to \mathbb{R}_+$ is defined by
		$$u^* (x_1,x_2,y) = \phi((|x|^{2\alpha+2} + (\alpha+1)^2y^2)^{1/2}) = \phi(r),$$
		where $r:= (|x|^{2\alpha+2} + (\alpha+1)^2y^2))^{1/2}$, $\phi: \mathbb{R}_+ \to \mathbb{R}_+$ and
		$$ \left| \{u^* > t\}\right|_{2,\alpha} = \left|\{u > t\}\right|_{2,\alpha}, \textrm{ for all } t>0. $$
	\end{definition}
	\begin{remark}\label{rem2}
		Put $M = \max_{\mathbb{R}^3} u$. It is not difficult to see that the map
		$$\lambda : t \mapsto |\{u>t\}|_{2,\alpha}$$
		is nonincreasing and right-continuous. Therefore the map
		$$ t \mapsto |\{t < u \leq  M\} \cap \{\nabla u = 0\}|_{2,\alpha} $$
		is nonincreasing and the function $\phi$ is nonincreasing, right-continuous. Moreover, the set
		$$ \{t \in \mathbb{R} : \exists s \in \mathbb{R}, \phi(s) = t, \phi'(s) = 0\} $$
		has Lebesgue measure 0 in $\mathbb{R}$. As in \cite{cia} the map $h : [0, M] \to [0, \infty)$ defined by
		$$h(t) = |\{t < u^* \leq M\} \cap \{\nabla u^* = 0\}|_{2,\alpha}$$
		is nonincreasing. Moreover $h'(t) = 0$ a.e. on $[0, M]$.
	\end{remark}
	
	We need the following lemma to prove Theorem \ref{thm2}.
	\begin{lemma}\label{lem1}
		Let $u \in C^\infty_0 (\mathbb{R}^3; \mathbb{R}_+)$. Assume that $u \not\equiv 0$. Denote $M = \max_{\mathbb{R}^3} u$. Then
		\begin{equation}\label{eqn8}
			\int\limits_{\{u^* =t\}} \dfrac{|x|^{2\alpha}}{|\nabla u^*|} d\mathcal{H}^2 \geq \int\limits_{\{u =t\}} \dfrac{|x|^{2\alpha}}{|\nabla u|} d\mathcal{H}^2 \textrm{ for a.e. } t\in [0,M].
		\end{equation}
	\end{lemma}
	
	\begin{proof}
		The proof is similar to that of \cite[Lemma 1]{nga} and for this reason, it is omitted here.
	\end{proof}
	
	Now we are able to give the proof of Theorem \ref{thm2}.
	
	\begin{proof}[Proof of Theorem \ref{thm2}]
		We follow the procedure of the proof of \cite[Theorem 2]{nga}. The case $u \equiv 0$ is simple. We may now assume $u \not\equiv 0$. Let $M = \max_{\mathbb{R}^3} u$. Notice that $u \in C^\infty_0 ( \mathbb{R}^3; \mathbb{R}_+)$. Then, by Sard's Theorem, the set
		$$\left\{ t\in[0,M]: \exists (x_1,x_2,y) \in \mathbb{R}^3 \textrm{ s.t. } u(x_1,x_2,y) = t, \textrm{ } \nabla u(x_1,x_2,y) = 0\right\}$$
		has Lebesgue measure 0 in $\mathbb{R}$. By the definition of $u^*$, the set
		$$\left\{ t\in[0,M]: \exists (x_1,x_2,y) \in \overline{\mathbb{R}}^3_{s_1} \textrm{ s.t } u^*(x_1,x_2,y) = t, \textrm{ } \nabla u^*(x_1,x_2,y) = 0\right\}$$
		has Lebesgue measure 0 in $\mathbb{R}$. Using the co-area formula, we deduce that
		$$ \int\limits_{\mathbb{R}^3} |\nabla_G u |^2 dx_1dx_2dy = \int\limits_0^M dt\int\limits_{u^{-1}(t)} |\nabla_G u |d\mu_G,$$
		and
		$$ \int\limits_{\mathbb{R}^3_{s_1}} |\nabla_G u^* |^2 dx_1dx_2dy = \int\limits_0^M dt\int\limits_{u^{*-1}(t)} |\nabla_G u^* |d\mu^*_G,$$
		where $d\mu_G = \frac{|\nabla_G u|}{|\nabla u|}d\mathcal{H}^2$ and $d\mu^*_G = \frac{|\nabla_G u^*|}{|\nabla u^*|}d\mathcal{H}^2$. To prove Theorem \ref{thm2}, we only need to show that for $t\in[0,M]$ such that $t$ is not a critical value of $u$ and $u^*$, the following inequality holds
		\begin{equation}\label{eqn9}
			\int\limits_{u^{*-1}(t)} |\nabla_G u^* |d\mu^*_G \leq \int\limits_{u^{-1}(t)} |\nabla_G u |d\mu_G.
		\end{equation}
		
		Notice that
		$$d\mu_G = \dfrac{|\nabla_G u|}{|\nabla u|}d\mathcal{H}^2 = \sqrt{\nu_1^2 + \nu_2^2 + |x|^{2\alpha} \nu_3^2} d\mathcal{H}^2.$$
		From this and the H\"older inequality, we obtain
		\begin{align}\label{eqn10}
			\nonumber
			\left(P_{2,\alpha} (\{u>t\})\right)^2 & = \Bigg(\int\limits_{u^{-1} \{t\}} |x|^\alpha d\mu_G \Bigg)^2 \\
			& \leq \Bigg(\int\limits_{u^{-1} \{t\}} |\nabla_G u| d\mu_G \Bigg)\Bigg(\int\limits_{u^{-1} \{t\}} \dfrac{|x|^{2\alpha}}{|\nabla_G u|} d\mu_G \Bigg).
		\end{align}
		By definition of $u^*$, we have
		$$ |\nabla_G u^* (x_1,x_2,y)| = (\alpha+1)|x|^\alpha |g'(r)|. $$
		Since $d\mu^*_G = \frac{|\nabla_G u^*|}{|\nabla u^*|}d\mathcal{H}^2$ and $r$ is constant on $u^{*-1} \{t\} \cap \mathbb{R}^3_{s_1}$, we deduce that
		\begin{align}\label{eqn11}
			\nonumber
			\left(P_{2,\alpha,1} (\{u^*>t\})\right)^2 & = \Bigg(\int\limits_{u^{*-1} \{t\}} |x|^\alpha d\mu^*_G \Bigg)^2 \\
			& = \Bigg(\int\limits_{u^{*-1} \{t\}} |\nabla_G u^*| d\mu^*_G \Bigg)\Bigg(\int\limits_{u^{*-1} \{t\}} \dfrac{|x|^{2\alpha}}{|\nabla_G u^*|} d\mu^*_G \Bigg).
		\end{align}
		Recall that 
		$$ |\{u>t\}|_{2,\alpha} = |\{u^* > t\}|_{2,\alpha}.$$
		Thus, Theorem \ref{thm1} implies that 
		\begin{equation}\label{eqn12}
			P_{2,\alpha,1} (\{u^*>t\}) \leq P_{2,\alpha} (\{u > t\}).
		\end{equation}
		The inequality \eqref{eqn9} then follows from Lemma \ref{lem1}, \eqref{eqn10}, \eqref{eqn11}, \eqref{eqn12} and the fact that $\frac{d\mu_G}{|\nabla_G u|} = \frac{d\mathcal{H}^2}{|\nabla u|}$. Our proof is complete. 
	\end{proof}
	
	\section{Sobolev-type inequalities with weights}\label{sec3}
	
	In this section, we will give the proof of Theorem \ref{thm3}. First, we need the following definition of the weighted Sobolev spaces.
	
	\begin{definition}\label{def5}
		Let $q>1$. We define $W_0^{1,2\alpha,q} (\mathbb{R}^3)$ as the completion of $C^\infty_0(\mathbb{R}^3)$ with respect to the norm
		$$\|u\|_{W_0^{1,2\alpha,q}} =  \Bigg( \int\limits_{\mathbb{R}^3} |\nabla_G u|^2dx_1dx_2dy \Bigg)^{1/2} + \Bigg(\int\limits_{\mathbb{R}^3} |x|^{2\alpha}|u|^q dx_1dx_2dy \Bigg)^{1/q}.$$
	\end{definition}
	
	For $u \in W_0^{1,2\alpha,q} (\mathbb{R}^3)\setminus \{0\}$, we consider the ratio
	\begin{equation}\label{eqn13}
		C_{2\alpha,q}(u) = \dfrac{\Bigg(\int\limits_{\mathbb{R}^3}|\nabla_G u|^2 dx_1dx_2dy\Bigg)^{1/2}}{\Bigg(\int\limits_{\mathbb{R}^3} |x|^{2\alpha}|u|^q dx_1dx_2dy \Bigg)^{1/q}}.
	\end{equation}
	By rescaling $X_1= \lambda x_1, X_2 = \lambda x_2$, and $Y = \lambda^{\alpha+1} y$, we have $U(x,y) = u(\lambda x, \lambda^{\alpha+1} y)$ and
	$$ \int\limits_{\mathbb{R}^3} |\nabla_G U|^2 dx_1dx_2dy = \lambda^{-(\alpha+1)}\int\limits_{\mathbb{R}^3} |\nabla_G u|^2 dX_1dX_2dY,$$
	and
	$$ \int\limits_{\mathbb{R}^3} |x|^{2\alpha}|U|^q dx_1dx_2dy = \lambda^{-(3\alpha+3)}\int\limits_{\mathbb{R}^3} |X|^{2\alpha} |u|^q dX_1dX_2dY,$$
	so that 
	$$ C_{2\alpha,q} (U) = \lambda^{\frac{3\alpha+3}{q} -\frac{\alpha+1}{2}} C_{2\alpha,q} (u).$$
	Hence, to have 
	$$ \inf\limits_{W_0^{1,2\alpha,q}(\mathbb{R}^3) \setminus \{0\} } C_{2\alpha,q} (u) > 0$$
	we need $q = 6$ to be the critical exponent.
	
	To prove Theorem \ref{thm3}, we need to recall the results of \cite[Lemma 2]{Tal}.
	\begin{lemma}\label{lem2}
		Let $m, p, q$ be real numbers such that
		$$1< p < m, \textrm{ } q = mp/(m-p).$$
		Let $\phi : \mathbb{R}_+ \to \mathbb{R}_+$ be a Lipschitz function and such that
		\begin{equation*}
			\int\limits_0^\infty r^{m-1} |\phi'(r)|^p dr < \infty, \textrm{ } \phi(r) \to 0 \textrm{ when } r \to \infty.
		\end{equation*}
		Then
		\begin{equation}\label{eqn14} \dfrac{\big(\int_0^\infty r^{m-1} |\phi'(r)|^p\big)^{1/p}}{\big(\int_0^\infty r^{m-1} |\phi(r)|^q\big)^{1/q}} \geq D_{p,q,m}\end{equation}
		with the best constant given by
		$$D_{p,q,m} = m^{\frac{1}{p}}\left(\dfrac{p-1}{m-1}\right)^{-\frac{1}{p'}} \left[ \dfrac{1}{p'} B\left( \dfrac{m}{p}, \dfrac{m}{p'} \right)\right]^{\frac{1}{m}},$$
		where $p' = p/(p-1)$. The equality sign in \eqref{eqn14} holds with
		$$ \phi (r) = (a + br^{p'})^{1-m/p},$$
		where $a, b$ are positive constants.
	\end{lemma}
	
	We also need some calculations. Using the polar coordinates
	$$x_1 = (r\sin\varphi)^{1/(\alpha+1)}\sin\theta, \textrm{ } x_2 = (r\sin\varphi)^{1/(\alpha+1)}\cos\theta, \textrm{ } y = \dfrac{r\cos \varphi}{\alpha+1},$$
	we have $dx_1dx_2dy = \frac{r^{2/(\alpha+1)}}{(\alpha+1)^2}(\sin \varphi)^{(1-\alpha)/(\alpha+1)}drd\theta d\varphi$. Next, consider the case $u(x_1,x_2,y) = \phi(r)$, where $r$ is given in Definition \ref{def4}. We have
	\begin{equation}\label{eqn15}
		\int_{\mathbb{R}^3_{s_1}} |x|^{2\alpha} |u|^{6}dx_1dx_2dy = \dfrac{2\pi}{n(\alpha)(\alpha+1)^2}  \int_0^\infty r^2|\phi(r)|^{6}dr,
	\end{equation}
	and
	\begin{equation}\label{eqn16}
		\int_{\mathbb{R}^3_{s_1}} |\nabla_G u|^2 dx_1dx_2dy =  \dfrac{2\pi}{n(\alpha)} \int_0^\infty r^2|\phi'(r)|^2dr.
	\end{equation}
	
	Now we are in a position to give the proof of Theorem \ref{thm3}.
	
	\begin{proof}[Proof of Theorem \ref{thm3}]
		The method of the proof is inspired by the one for the two-dimensional case in \cite{nga}, which makes use of the P\'olya-Szeg\"o type inequality. Since $C^\infty_0(\mathbb{R}^3)$ is dense in $W^{1,2\alpha,6}_0(\mathbb{R}^3)$, we can assume that $u \in C^\infty_0(\mathbb{R}^3)$. Notice that 
		$$|\nabla_G |u|| \leq |\nabla_G u|$$
		thus we can assume $u$ is a nonnegative function. It is easy to see that we can approximate a nonnegative function $w \in W^{1,2\alpha,6}_0(\mathbb{R}^3)$ by a sequence of nonnegative functions $\{w_j\}_{j=1}^\infty$, we only need to give the proof for $u \in C^\infty_0(\mathbb{R}^3)$ and $u\geq 0$.
		
		Let $u^*$ be the rearrangement of $u$. The property of rearrangement infers that
		\begin{align}\label{eqn17}
			\int_{\mathbb{R}^3_{s_1}} |x|^{2\alpha} |u^*|^{6} dx_1dx_2dy = \int_{\mathbb{R}^3} |x|^{2\alpha} |u|^{6} dx_1dx_2dy.
		\end{align}
		Note that $u \in C^\infty_0(\mathbb{R}^3)$, $u\geq 0$. It follows from Theorem \ref{thm2} that
		\begin{equation}\label{eqn18}
			\int_{\mathbb{R}^3_{s_1}} |\nabla_G u^*|^2 dx_1dx_2dy \leq \int_{\mathbb{R}^3} |\nabla_G u|^2dx_1dx_2dy.
		\end{equation}
		
		On the other hand, we have
		$$u^*(x_1,x_2,y) = \phi(r).$$ 
		Note that the function $\phi$ satisfies Lemma \ref{lem2}. Then, from \eqref{eqn15}, \eqref{eqn16} and Lemma \ref{lem2}, we deduce that
		\begin{equation}\label{eqn19}
			L_{2\alpha,6} \bigg( \int_{\mathbb{R}^3_{s_1}} |x|^{2\alpha} |u^*|^{6} dx_1dx_2dy\bigg)^{1/6} \leq \bigg( \int_{\mathbb{R}^3_{s_1}} |\nabla_G u^*|^2dx_1dx_2dy\bigg)^{1/2},
		\end{equation}
		where
		\begin{equation}\label{eqn20}
			L_{2\alpha,6} = \left( \dfrac{2\pi}{n(\alpha)}\right)^{\frac{-1}{3}} (\alpha+1)^{ \frac{-1}{3}} D_{2,6,3} 
		\end{equation}
		with $D_{2,6,3}$ is the constant given in Lemma \ref{lem2}. From \eqref{eqn17} - \eqref{eqn20}, our Sobolev-type inequality \eqref{eqn3} follows. Our proof is complete.
	\end{proof}
	
	\begin{remark}
		\label{rem3}
		Our proof only gives a lower bound for the best constant $C_{2\alpha,6}$. The exact value of $C_{2\alpha,6}$ remains an open question.
	\end{remark}

	\section{Applications to a class of degenerate elliptic equations}\label{sec4}
	In this section, we are mainly concerned with the existence and nonexistence of nontrivial solutions to the problem \eqref{P}. Notice that our methods are similar to the ones that have been used in \cite{Luyen2}, and thus many details will be omitted. 
	
	\subsection{Nonexistence result}\label{sec41}
	In this subsection, we will derive the nonexistence of nontrivial solutions to the problem \eqref{P} for $f(x_1,x_2, y, \xi) = |x|^{2\alpha} |\xi|^{p-1} \xi$, $p  \geq 1$. To this end, let
	$$ F(x_1,x_2,y,\xi) := \int\limits_0^\xi f(x_1,x_2,y,\tau)d\tau = \dfrac{|x|^{2\alpha}}{p+1} |\xi|^{p+1}. $$
	
	\begin{definition}\label{defextra}
		We denote by $\mathcal{S}_2(\overline{\Omega})$ the linear space of functions in $C_0^1 (\overline{\Omega})$ such that
		$$\dfrac{\partial^2 u}{\partial x_1^2}, \textrm{ } \dfrac{\partial^2 u}{\partial x_2^2}, \textrm{ and } |x|^{2\alpha} \dfrac{\partial^2 u}{\partial y^2} \textrm{ (in distribution sense)}$$ 
		are continuous in $\Omega$ and can be continuously extended to $\overline{\Omega}$.
	\end{definition}
	
	A function $u(x_1,x_2,y) \in \mathcal{S}_2(\overline{\Omega})$ is said to be a solution to the problem \eqref{P} if
	\begin{align*}
		\begin{cases}
			- \Delta_x u - |x|^{2\alpha} \dfrac{\partial^2 u}{\partial y^2} = |x|^{2\alpha}|u|^{p-1}u & \textrm{in } \Omega, \\
			u = 0 & \textrm{on } \partial \Omega.
		\end{cases}
	\end{align*}

	\begin{lemma}\label{lem3}
		Let $u(x_1,x_2,y)\in\mathcal{S}_2(\overline{\Omega})$ be a solution of \eqref{P}. Denote by $\nu = (\nu_1,\nu_2,\nu_3)$ the unit outward normal on $\partial\Omega$. Then we have
		\begin{align}
			\label{eqn21}
			\nonumber
			\bigg(\dfrac{3\alpha+3}{p+1}-&\dfrac{\alpha+1}{2}\bigg)\int_\Omega |x|^{2\alpha}|u|^{p+1}dx_1dx_2dy \\
			& = \int_{\partial\Omega} [x_1\nu_1+x_2\nu_2 + (1+\alpha)y\nu_3](\nu_1^2+\nu_2^2 + \nu_3^2 |x|^{2\alpha}) \bigg( \dfrac{\partial u }{\partial \nu}\bigg)^2 ds.
		\end{align}
	\end{lemma}
	
	\begin{proof}
		The proof is similar to that of the lemma in \cite{Tri3}, which is based on a series of straight forward calculation involving integration by parts, and for this reason, it is omitted here.
	\end{proof}
	
	\begin{definition}\label{def6}
		A domain $\Omega$ is called $G_\alpha$ -- star-shaped with respect to the origin if $(0, 0, 0) \in \Omega$ and $x_1\nu_1 +x_2\nu_2 + (1+\alpha)y\nu_3 \geq 0$ at every point of $\partial\Omega$.
	\end{definition}
	
	\begin{proof}[Proof of Theorem \ref{thm4}]
		The desired result is obtained directly from Lemma \ref{lem3} and Definition \ref{def6}.
	\end{proof}

	\subsection{Existence result}\label{sec42}
	In this subsection, we will give the existence result for the problem \eqref{P}. 
	\begin{definition}\label{def7}
		Let $\Omega$ be a bounded domain in $\mathbb{R}^3$. We denote  by $L^p_{|x|^{2\alpha}} (\Omega)$ the set of all measurable functions $u:\Omega \to \mathbb{R}$ such that 
		$$\int\limits_\Omega |x|^{2\alpha} |u|^p dx_1dx_2dy < \infty. $$
		We define the norm in $L^p_{|x|^{2\alpha}} (\Omega)$ as follows
		$$\|u\|_{L^p_{|x|^{2\alpha}} (\Omega)} := \bigg(\int\limits_\Omega |x|^{2\alpha} |u|^p dx_1dx_2dy\bigg)^{1/p}$$
		for every $u \in L^p_{|x|^{2\alpha}} (\Omega)$.
	\end{definition}
	
	\begin{definition}\label{def8}
		Let $\Omega$ be a bounded domain in $\mathbb{R}^3$. We denote by $\mathcal{S}^2_1 (\Omega)$ the set of all functions $u\in L^2(\Omega)$ such that 
		$$\dfrac{\partial u}{\partial x_1}, \textrm{ } \dfrac{\partial u}{\partial x_2}, \textrm{ and } |x|^\alpha\dfrac{\partial u}{\partial y} \in L^2(\Omega).$$
		We define the norm in $\mathcal{S}^2_1 (\Omega)$ as follows
		$$ \|u \|_{\mathcal{S}^2_1 (\Omega)} := \bigg( \int_\Omega |u|^2 dx_1dx_2dy \bigg)^{1/2} + \bigg(\int_\Omega |\nabla_G u|^2dx_1dx_2dy\bigg)^{1/2}. $$
		We also define the inner product in $\mathcal{S}^2_1 (\Omega)$ by
		$$(u,v)_{\mathcal{S}^2_1 (\Omega)} := (u,v)_{L^2(\Omega)} + (\nabla_G u, \nabla_G v)_{L^2(\Omega)}.$$
	\end{definition}
	
	\begin{definition}\label{def9}
		The space $\mathcal{S}^2_{1,0}(\Omega)$ is defined as the closure of $C^1_0 (\Omega)$ in the space $\mathcal{S}^2_1 (\Omega)$.
	\end{definition}
	
	\begin{remark}\label{rem4}
		By a similar argument as in the proof of \cite[Theorem 6]{Tri3}, we also have the two norms $\|u\|_{\mathcal{S}^2_1 (\Omega)}$ and
		$$\|u\|_{\mathcal{S}^2_{1,0}(\Omega)} := \bigg( \int_\Omega |\nabla_G u|^2 dx_1dx_2dy\bigg)^{1/2}$$
		are equivalent in $\mathcal{S}^2_{1,0}(\Omega)$.
	\end{remark}
	
	We need the following embedding theorem.	
	\begin{theorem}\label{thm6}
		Let $\Omega$ be a bounded domain in $\mathbb{R}^3$ with smooth boundary $\partial \Omega$ such that $(0,0,0) \in \Omega$. Then the embedding
		$$\mathcal{S}^2_{1,0} (\Omega) \hookrightarrow L^q_{|x|^{2\alpha}} (\Omega), \textrm{ where } 1\leq q \leq 6,$$
		is continuous, i.e., there exists a constant $C_q > 0$ such that 
		$$ \|u\|_{L^q_{|x|^{2\alpha}}(\Omega)} \leq C_q \|u\|_{\mathcal{S}^2_1(\Omega)}, \textrm{ } \forall u \in \mathcal{S}^2_{1,0} (\Omega).$$
		Moreover, the embedding
		$$\mathcal{S}^2_{1,0} (\Omega) \hookrightarrow L^q_{|x|^{2\alpha}} (\Omega), \textrm{ where } 1\leq q < 6,$$
		is compact.
	\end{theorem}
	
	\begin{proof}
		The proof of this lemma is similar to that of \cite[Theorem 3.2]{Luyen2}, which is based on Theorem \ref{thm3} (for the continuity of the embedding) and Fr\'echet-Kolmogorov Theorem (for the compactness). We omit the details.
	\end{proof}
	
	Now we are able to study the existence result for the problem \eqref{P}. 
	
	\begin{definition}\label{def10}
		A function $u \in \mathcal{S}^2_{1,0} (\Omega)$ is called a weak solution of the problem \eqref{P} if the following identity
		$$ \int_\Omega \nabla_G u \cdot \nabla_G \varphi dx_1dx_2dy - \int_\Omega f(x,y,u)\varphi dx_1dx_2dy = 0$$
		holds for every $\varphi \in \mathcal{S}^2_{1,0}(\Omega)$.
	\end{definition}
	
	The following lemma is obtained from the H\"older inequality and Theorem \ref{thm6}.	
	\begin{lemma}\label{lem4}
		Assume that $f : \Omega \times \mathbb{R} \to \mathbb{R}$ is a Carath\'eodory function such that there exist $q \in (2, 6), f_1(x, y) \in L^{p_1}_{|x|^{2\alpha}}(\Omega, \mathbb{R}_+), f_2(x, y) \in L^{p_2}_{|x|^{2\alpha}}(\Omega, \mathbb{R}_+)$, where $p_2 > 1$, $qp_2/(p_2-1) \leq 6$, and $p_1 > \frac{6p_2}{p_2(q-1)+6}$ such that
		$$|f(x_1,x_2, y, \xi)| \leq |x|^{2\alpha} (f_1(x_1,x_2, y) + f_2(x_1,x_2, y) |\xi|^{q-1})$$
		for a.e.  $(x_1,x_2, y) \in \Omega, \forall \xi \in \mathbb{R}$ .Then $\Psi_1(u) \in C^1(\mathcal{S}^2_{1,0}(\Omega), \mathbb{R})$ and 
		$$\Psi'_1(u)(v) = \int_\Omega f(x_1,x_2, y, u)vdx_1dx_2dy$$
		for all $v \in \mathcal{S}^2_{1,0}(\Omega)$, where
		$$\Psi_1 (u) = \int_\Omega F(x_1,x_2, y, u) dx_1dx_2dy,$$
		and $F(x_1,x_2, y, \xi) =\int_0^\xi f(x_1,x_2, y, \tau)d\tau$.
	\end{lemma}
	
	In the following we define the energy functional $\Phi: \mathcal{S}^2_{1,0}(\Omega) \to \mathbb{R}$ associated with the problem \eqref{P} by letting
	\begin{equation}\label{eqn22}
		\Phi (u) := \dfrac{1}{2}\int_\Omega |\nabla_G u|^2 dx_1dx_2dy - \int_\Omega F(x,y,u)dx_1dx_2dy
	\end{equation}
	for all $u \in \mathcal{S}^2_{1,0}(\Omega)$, where $F$ is defined as in Lemma \ref{lem4}. It follows from Lemma \ref{lem4} and the condition (A1) that $\Phi$ is well-defined on $\mathcal{S}^2_{1,0}(\Omega)$ and $\Phi \in C^1(\mathcal{S}^2_{1,0}(\Omega), \mathbb{R})$ with
	$$\langle\Phi'(u), v \rangle = \int_\Omega \nabla_G u \cdot \nabla_G v dx_1dx_2dy - \int_\Omega f(x,y,u)vdx_1dx_2dy $$
	for all $v \in \mathcal{S}^2_{1,0} (\Omega)$. Therefore, the weak solutions of the problem \eqref{P} are critical points of the functional $\Phi$.
	
	We need to use Mountain Pass Lemma to study the existence of critical points of the functional $\Phi$. First, let us recall the notion of the $(C)_c$ condition, which is an important notion in the statement of this lemma.
	
	\begin{definition}\label{def11}
		Let $\boldsymbol{X}$ be a real Banach space with its dual space $\boldsymbol{X}^*$ and let $\Psi \in C^1(\boldsymbol{X}, \mathbb{R})$. For $c \in \mathbb{R}$ we say that $\Psi$ satisfies the $(C)_c$ condition if for any sequence $\{u_n\}_{n=1}^\infty \subset \boldsymbol{X}$ with 
		$$\Psi(u_n) \to c \textrm{ and  }(1 + \|u_n\|_{\boldsymbol{X}}) \|\Psi'(u_n)\|_{\boldsymbol{X}^*} \to 0,$$
		there exists a subsequence $\{u_{n_k}\}_{k=1}^\infty$ that converges strongly in $\boldsymbol{X}$. 
	\end{definition}
	
	We will need the following version of Mountain Pass Lemma (see also \cite{cer1, cer2}).
	
	\begin{lemma}\label{lem5}
		Let $\boldsymbol{X}$ be a real Banach space and let $\Psi \in C^1(\boldsymbol{X}, \mathbb{R})$ satisfy the $(C)_c$ condition for any $c \in \mathbb{R}$, $\Psi(0) = 0$ and
		\begin{itemize}
			\item[(i)] There exist constants $\rho, \alpha > 0$ such that $\Psi(u) \geq \alpha, \forall u \in \boldsymbol{X}, \|u\|_{\boldsymbol{X}} = \rho$; 
			\item[(ii)] There exists an $u_1 \in \boldsymbol{X}$, $\|u_1\|_{\boldsymbol{X}} \geq \rho$ such that $\Psi(u_1) \leq 0$.
		\end{itemize}
		Then $\beta := \inf\limits_{\lambda\in\Lambda}\max\limits_{0\leq t\leq1}\Psi(\lambda(t)) \geq \alpha$ is a critical value of $\Psi$, where
		$$\Lambda := \{\lambda \in C ([0; 1], \boldsymbol{X}) : \lambda(0) = 0, \lambda(1) = u_1\}.$$
	\end{lemma}
	
	To conclude our paper, we have the following.
	
	\begin{proof}[Proof of Theorem \ref{thm5}]
		The proof follows the same arguments presented in the proof of \cite[Theorem 4.3]{Luyen2}, and for this reason, the details will be omitted here. The idea is that, based on the fact that $f$ satisfies the conditions (A1)-(A5) and Theorem \ref{thm6}, one can verify that the functional $\Phi$ satisfies the conditions of Lemma \ref{lem5}. Thus, our proof is complete. 
	\end{proof}
	
	\section*{Data availability}
	Not applicable.
	
	\section*{Declarations}
	The authors have no conflicts of interest to declare that are relevant to the content of this article.
	
	\bibliographystyle{acm}

\end{document}